\documentclass[11pt,a4paper]{amsart}%
\usepackage{amsfonts,amsmath,amssymb}
\usepackage{hyperref}
\usepackage{amssymb,amsthm,amsxtra}
\usepackage[usenames]{color}
\usepackage{amscd}
\usepackage{amsthm}
\usepackage{amsfonts}
\usepackage{amssymb}
\usepackage{amsmath}
\usepackage{graphicx}
\usepackage{enumerate}%
\usepackage[all]{xy}
\usepackage{aurical}
\usepackage[T1]{fontenc}
\usepackage{multicol}
\let\oldtocsubsection=\tocsubsection
\let\oldtocsubsubsection=\tocsubsubsection
\renewcommand{\tocsubsection}[2]{\hspace{1em}\oldtocsubsection{#1}{#2}}
\renewcommand{\tocsubsubsection}[2]{\hspace{2em}\oldtocsubsubsection{#1}{#2}}

\newcommand*\IsoTo{%
  \xrightarrow[]{\raisebox{-0.5 em}{\smash{\ensuremath{\sim}}}}%
}

\theoremstyle{plain}

\newtheorem*{theorem*}{Theorem}
\newtheorem{lemma}{Lemma}[section]

\newtheorem*{conjecture*}{Conjecture}

\newtheorem{thm}[lemma]{Theorem}

\newtheorem{lem}[lemma]{Lemma}
\newtheorem{cor}[lemma]{Corollary}

\newtheorem{introtheorem}{Theorem}

\newtheorem{introcor}[introtheorem]{Corollary}
\newtheorem{introthm}[introtheorem]{Theorem}

\theoremstyle{definition}

\newtheorem{definition}[lemma]{Definition}
\newtheorem{defn}[lemma]{Definition}

\newtheorem{notn}[lemma]{Notation}

\newtheorem{example}[lemma]{Example}

\newtheorem*{example*}{Example}

\theoremstyle{remark}

\newtheorem*{remark*}{Remark}

\newtheorem{remark}[lemma]{Remark}

 \theoremstyle{plain}

\newcommand{\tr}{\operatorname{Tr}}

\newcommand{\ind}{\operatorname{ind}}

\newcommand{\Rep}{\operatorname{Rep}}

\newcommand{\oH}{\operatorname{H}}

\newcommand{\diag}{\operatorname{diag}}

\newcommand{\A}{\mathbb{A}}

\newcommand{\eps}{\varepsilon}

\newcommand{\Ker}{\operatorname{Ker}}

\newcommand{\Z}{{\mathbb Z}}
\newcommand{\Q}{{\mathbb Q}}
\newcommand{\R}{{\mathbb R}}

\newcommand{\C}{{\mathbb C}}

\newcommand{\Exp}{\operatorname{Exp}}

\newcommand{\Id}{\operatorname{Id}}

\newcommand{\Fre}{{Fr\'{e}chet \,}}

\newcommand{\bfG}{{\mathbf{G}}}

\newcommand{\bfN}{{\mathbf{N}}}
\newcommand{\bfU}{{\mathbf{U}}}

\newcommand{\cF}{{\mathcal{F}}}

\newcommand{\cN}{{\mathcal{N}}}

\newcommand{\g}{{\mathfrak{g}}}

\newcommand{\fg}{{\mathfrak{g}}}
\newcommand{\fh}{{\mathfrak{h}}}
\newcommand{\fn}{{\mathfrak{n}}}
\newcommand{\fk}{{\mathfrak{k}}}
\newcommand{\fu}{{\mathfrak{u}}}
\newcommand{\fv}{{\mathfrak{v}}}
\newcommand{\fw}{{\mathfrak{w}}}

\newcommand{\cO}{{\mathcal{O}}}
\newcommand{\GL}{\operatorname{GL}}

\newcommand{\SL}{\operatorname{SL}}
\newcommand{\Sp}{\operatorname{Sp}}

\newcommand{\SO}{\operatorname{SO}}

\newcommand{\cM}{\mathcal{M}}

\newcommand{\fp}{\mathfrak{p}}
\newcommand{\fl}{\mathfrak{l}}
\newcommand{\fm}{\mathfrak{m}}

\newcommand{\fr}{\mathfrak{r}}

\newcommand{\cV}{\mathcal{V}}
\newcommand{\cW}{\mathcal{W}}

\newcommand{\Dima}[1]{{{}{#1}}}
\newcommand{\DimaA}[1]{{{#1}}}
\newcommand{\DimaB}[1]{{{#1}}}

\newcommand{\DimaC}[1]{{{#1}}}
\newcommand{\DimaD}[1]{{{#1}}}

\newcommand{\onto}{{\twoheadrightarrow}}

\newcommand{\into}{{\hookrightarrow}}

\newcommand{\ot}{\leftarrow}

\newcommand{\Orb}{\mathcal{O}}
\newcommand{\WF}{\operatorname{WF}}
\newcommand{\WO}{\operatorname{WO}}
\newcommand{\WS}{\operatorname{WS}}

\renewcommand{\sl}{\mathfrak{sl}}
\newcommand{\F}{{F}}

\begin{document}


\author{Dmitry Gourevitch}
\address{ Department of Mathematics,
Weizmann Institute of Science,
 Rehovot 7610001 Israel}
\email{dmitry.gourevitch@weizmann.ac.il}

\author{Siddhartha Sahi}
\address{Department of Mathematics, Rutgers University, Hill Center -
Busch Campus, 110 Frelinghuysen Road Piscataway, NJ 08854-8019, USA}
\email{sahi@math.rugers.edu}

\date{\today}
\title{Generalized and degenerate Whittaker quotients and Fourier coefficients} 


\begin{abstract}

The study of Whittaker models for representations of reductive groups over local and global fields has become a central tool in representation theory and the theory of automorphic forms. However, only generic representations have Whittaker models. In order to encompass other representations, one attaches a degenerate (or a generalized) Whittaker model $\cW_{\cO}$, or a Fourier coefficient in the global case, to any nilpotent orbit.

In this note we survey some classical and some recent work in this direction - for Archimedean, p-adic and global fields. The main results concern the existence  of models. For a representation $\pi$, call the set of maximal \DimaD{elements of the set of} orbits $\cO$ with $\cW_{\cO}$ that includes $\pi$ the \emph{Whittaker support} of $\pi$. The two main questions discussed in this note are:
\begin{enumerate}
\item What kind of orbits can appear in the Whittaker support of a representation?
\item How does the Whittaker support of a given representation $\pi$ relate to other invariants of $\pi$, such as its wave-front set?
\end{enumerate}

%
%
%
%
\end{abstract}

\maketitle
\tableofcontents


%

\section{Introduction}

The study of Whittaker and degenerate Whittaker models for
representations of reductive groups over local fields evolved in
connection with the theory of automorphic forms (via their Fourier
coefficients), and has found important applications in both areas.
See for example \cite{GK,Sh74,NPS73,Kos,CHM} for Whittaker models and \cite{Ka85,Ya86, WaJI,MW,Mat,Ginz,Ji07,GRS_Book,Ginz2,GZ,JLS,CFGK} for degenerate and generalized Whittaker models and Fourier coefficients.

In this note we survey some classical results and some recent work on
degenerate and generalized Whittaker models and Fourier coefficients, including some recent results of the authors  obtained in collaboration with R. Gomez and others.

Let $F$ be  either $\R$ or a finite extension of $\Q_p$, and let $G$ be a finite central extension of the group of $F$-points of a connected reductive algebraic group defined over $F$\footnote{We view complex reductive groups as a special case of real reductive groups.}.
Let $\Rep^{\infty}(G)$ denote the category of smooth representations of $G$ (see \S \ref{subsec:not} below).
Let $\fg$ denote the Lie algebra of $G$ and $\fg^*$ denote its dual space. To every coadjoint nilpotent orbit $\cO\subset \fg^*$ and every $\pi\in\Rep^{\infty}(G)$ we associate a certain generalized Whittaker quotient $\pi_{\cO}$ (see \S \ref{subsec:DefMod} below). Let $\WO(\pi)$ denote the set of all nilpotent orbits $\cO$ with $\pi_{\cO}\neq 0$ and $\WS(\pi)$ denote the set of maximal orbits in $\WO(\pi)$  with respect to the closure ordering. We call $\WS(\pi)$ the \emph{Whittaker support of $\pi$.}

Denote by $\cM(G)$ the category of admissible (finitely-generated) representations (see \S \ref{subsec:not} below). For $\pi\in \cM(G)$, and a nilpotent orbit $\cO\subset \fg^*$, \cite{HowGL,HCWF,BV} define a coefficient $c_{\cO}(\pi)$ using the asymptotics of the character of $\pi$ at $1\in G$ (see \S \ref{sec:WF} below). Denote by $\WF(\pi)$ the set of maximal elements in the set of orbits with non-zero coefficients. For non-Archimedean $F$ this  set coincides with $\WS(\pi)$.

\begin{introthm}[{\cite[Proposition I.11, Theorem I.16 and Corollary I.17]{MW}, and \cite{Var}}] \label{thm:MW}
Assume that $F$ is non-Archimedean and  $G$ is algebraic\footnote{It seems that the assumption that $G$ is algebraic is not used in the proof.} and let $\pi\in \cM(G)$. Then
\begin{enumerate}[(i)]
\item $\WF(\pi)=\WS(\pi)$.
\item For any $\cO\in \WF(\pi), \, c_{\cO}(\pi)=\dim \pi_{\cO}$.
\end{enumerate}
\end{introthm}

We conjecture that the same  holds for Archimedean $F$. In \S \ref{sec:WF} we define the wave-front set $\WF(\pi)$, sketch the proof of Theorem \ref{thm:MW}  and
survey some partial results in the Archimedean case, based on \cite{GW,Mat,MatENS,GGS,GGS2,GSS}.

The next question that arises is what orbits can appear in Whittaker supports of representations. Based on the Kirillov orbit method one may conjecture that they are all admissible. This notion has to do with splitting of a certain metaplectic double cover of  the centralizer $G_{\varphi}$ for any $\varphi$ in the orbit (see \S \ref{subsec:cov} below).
Based on the Langlands  correspondence it makes sense to conjecture that if $G$ is algebraic then these orbits are special in the sense of Lusztig (\cite{LusSpec}, \cite[Ch. III]{Spa}). More precisely, an orbit is said to be special if the corresponding orbit over the algebraic closure is special.
For non-Archimedean $F$ these conjectures might hold, but for Archimedean $F$ they do not. For example, the minimal orbit of $G_2(\R)$ is not special and appears in $\WS(\pi)$ for  some irreducible unitary $\pi$ by \cite{VogG2}.
Also, the minimal orbit of $U(2,1)$ is non-admissible and appears in $\WS(\pi)$ for some $\pi$.
 On the positive side, the following holds for all $F$.

\begin{introthm}[{\cite[Theorem A]{GGS2}, following \cite{Mog,GRS,Ginz,JLS}}]\label{thm:adm}
Let $\pi\in \Rep^{\infty}(G)$ and let $\cO\in \WS(\pi)$. Then
$\cO$ is a quasi-admissible orbit (see \S \ref{subsec:cov} below).
\end{introthm}

Also, for non-Archimedean $F$ and any $\pi \in \Rep^{\infty}(G_2(F))$  all $\cO\in \WS(\pi)$ are special by \cite{JLS,LS}.

Let us summarize the relations between the above notions for classical groups.
\begin{introthm}[{\cite[\S 6]{GGS2}, based on \cite{Mog,Nev,Oht}}]\label{prop:ClasQuas}
Let $\cO\subset \fg^*$ be a nilpotent orbit.
\begin{enumerate}[(i)]
\item If $G$ is $\GL_n(F)$ or $\SL_n(F)$ then all orbits are admissible and special.
\item If $G$ is $U(V)$ or $SU(V)$ for a hermitian space $V$ then all orbits are quasi-admissible and special. If $F\neq \R$ then they are  all admissible.
\item If $G$ if either $\Sp_{2n}(F)$, or $O(V)$ or $\SO(V)$ (for a quadratic space $V$ over $F$),  then the following are equivalent:
$$ (a) \,\,  \cO \text{ is admissible}\quad (b)\,\, \cO \text{ is quasi-admissible}  \quad (c)\,\, \cO \text{ is special} $$
\item If \DimaD{$F=\R$}  and $G$ is a classical group not listed above then all orbits are admissible.
\end{enumerate}
\end{introthm}

It is possible  that the notions of admissible and quasi-admissible are equivalent for all $G$ in the case when $F$ is non-archimedean. These notions differ for $U(p,q)$ and for $SU(p,q)$. 
They also differ for the split real forms of $E_7$ and $E_8$, though we do not know whether the non-admissible quasi-admissible orbits appear in Whittaker supports of representations. It is also possible that all special orbits are quasi-admissible for all groups. If $G$ is algebraic, $F$ is Archimedean and $\pi\in \cM(G)$ has integral infinitesimal character then all  $\cO\in \WF(\pi)$ are special, cf. \cite[Theorem D]{BVClass} and \cite[Theorem 1.1]{BVExc}.

Under additional assumptions on $\pi$ one can show that all the  \DimaD{orbits in $\WS(\pi)$} are $F$-distinguished, i.e. do not intersect the Lie algebras of proper Levi $F$-subgroups of $G$.

\begin{introthm}[{\cite{Mog,Harris,GGS2}}] \label{thm:cuspTemp}
Let $\pi\in \Rep^{\infty}(G)$,
 and let $\cO\subset \fg^*$ be a nilpotent orbit.
 \begin{enumerate}[(i)]
 \item \label{it:cusp} If $F$ is non-Archimedean, $\cO \in \WS(\pi)$ and $\pi$ is quasi-cuspidal then $\cO$ is $F$-distinguished.
\item If $\pi$ is admissible and tempered, $\cO\in\WF(\pi)$ and either $G$ is classical or $F$ is Archimedean then $\cO$ is $F$-distinguished.
\end{enumerate}
\end{introthm}

Note that if $G$ is semi-simple then $\cO$ is $F$-distinguished if and only if all reductive subgroups of the centralizer of any element of $\cO$ are compact. Thus, such orbits are sometimes called compact. In the Archimedean case, compact orbits were classified in \cite{PT}.

Over Archimedean fields, all the orbits in $\WF(\pi)$ lie in the same complex orbit \Dima{by \cite[Theorem D]{Ross}.} This is conjectured to hold for p-adic $F$ as well, but so far \Dima{this conjecture} is known only for the group $\GL_n(F)$ (see \cite[\S II.2]{MW}).
Over finite fields, an analogous result is proven in \cite{LusFin}, as well as some  analogues of Theorems \ref{thm:MW} and \ref{thm:adm}.
%
%

In \cite{GZ,LM,Prz} it is shown that the Whittaker support and the  wave-front set have the expected behaviour  under $\theta$-correspondence.

The behaviour of wave-front set under parabolic and cohomological induction  is studied in \cite{MW,BarBoz,BVInd} (see \S \ref{subsec:DefWF} below).

In the non-Archimedean case, the functor $\pi \mapsto \pi_{\cO}$ is exact. We conjecture that the same holds in the Archimedean case, at least on the subcategory of admissible representations  $\pi$ with $\WO(\pi)$ lying in the closure of $\cO$. This is shown for principal nilpotent $\cO$ in \cite{CHM}. For some partial results in this direction for other $\cO$ see \cite{L, WaJI,MatENS,WaJI,AGS2,GGS}.

One can define more general  models and quotients. Namely, let $S\in \fg$ such that the adjoint action of $S$ is diagonalizable over $\Q$, and let $\varphi\in \fg^*$ such that $ad^*(S)(\varphi)=-2\varphi$. In \S \ref{subsec:DefMod} below we define a \emph{degenerate Whittaker model} $\cW_{S,\varphi}$, and define a degenerate Whittaker quotient of $\pi$ to be the coinvariants $\pi_{S,\varphi}:=(\cW_{S,\varphi}\otimes \pi)_G$. In \S \ref{subsec:compar} we discuss an epimorphism $\cW_{\cO}\onto \cW_{S,\varphi}$ constructed in \cite{GGS}, under the condition that the orbit $\cO$ intersects the closure of the orbit of $\varphi$ under the centralizer of $S$ in $G$. The existence of this epimorphism sheds some light on the non-maximal orbits in $\WO(\pi)$. In particular, it is shown that for $\pi\in \Rep^{\infty}(\GL_n(F))$, we have $\cO\in \WO(\pi)$ if and only if there exists an orbit $\cO'\in \WS(\pi)$ with $\cO\subset \overline{\cO'}$.
\DimaA{The analogous statement for other groups is not known in general}.
For  $\pi\in \cM(\GL_n(F))$ we also have $\WS(\pi)=\WF(\pi)$.
Also,  for $\cO\in \WS(\pi),$ the induced map $\pi_{\cO}\onto \pi_{S,\varphi}$ is non-zero and is an isomorphism for non-archimedean $F$. This fact is used in the proof of Theorems \ref{thm:adm} and \ref{thm:cuspTemp}\eqref{it:cusp}.

Over the adeles one defines the Whittaker support using period integrals rather than quotients \Dima{(see \S \ref{subsec:FC} below for more details). These period integrals are called \emph{Fourier coefficients}. We refer the reader to \cite{Ginz} for some intriguing open questions on Fourier coefficients of automorphic forms.  Most of the questions described above work analogously over the adeles.}
Let $K$ be a global field, $\mathbf G$ be a reductive group defined over $K$, and $\pi$ be an automorphic representation of the adelic points $G:=\mathbf{G}(\A)$. 

\begin{introthm}[{\cite{GRS_Sp,Ginz,JLS,Shen,GGS2}}]\label{thm:Glob}
 Let  $\cO \in \WS(\pi)$ be a nilpotent $\mathbf{G}(K)$-orbit. Then
\begin{enumerate}[(i)]
\item \label{it:GlobDeg} $\cO$ is quasi-admissible. Furthermore, if  $\bfG$ is classical then  $\cO$ is special.
\item \label{it:GlobCusp} If $\pi$ is cuspidal then $\varphi$ does not belong to the Lie algebra of any Levi subgroup of $\bfG(K)$ defined over $K$.
\end{enumerate}
\end{introthm}

\begin{introcor}\label{cor:GlobCusp}
 Let  $\cO \in \WS(\pi)$ be a nilpotent $\mathbf{G}(K)$-orbit. Assume that $\pi$ is cuspidal and fix  $\varphi \in  \cO$. Then any Levi subgroup of the stabilizer of $\varphi$ in $\bfG(K)$ is $K$-anisotropic. Moreover, assume that $\bfG$ is classical, and let $\lambda$ be the  partition corresponding to $\varphi$. Then
\begin{enumerate}[(i)]
\item If $\bfG=\GL_n$ or $\bfG=\SL_n$ then $\lambda$ consists of one part, i.e. $\pi$ is generic.
\item  If $\bfG=\Sp_n$  then $\lambda$ is totally even.
\item  If $\bfG=\SO_n$ or $\bfG=\mathrm{O}_n$ then $\lambda$ is totally odd, i.e. consists of odd parts only.
\end{enumerate}
\end{introcor}

 The study of Fourier coefficients of automorphic representations has applications in string theory, cf.  \cite{GMV} and \cite[Ch. 2]{FGKP}. A special role is played by Fourier coefficients of minimal and next-to-minimal representations \DimaB{of the split simply laced groups $\SL(n),SO(5,5), E_6,E_7,E_8$. In  \cite{GSPhys2} we  express any minimal \DimaD{or} next-to-minimal automorphic form for these groups through its Whittaker-Fourier coefficients, \emph{i.e.} period integrals over the nilradical of a Borel subgroup of $G$ against a character of this subgroup, following \cite{MilSah,AGKLP}. This  is important since the  latter integral is \DimaD{often} Eulerian\DimaD{, which allows it to be computed explicitly.}. On the other hand, the analogous statement does not hold for $Sp(4)$. We do not expect it to hold for higher symplectic groups either.

\Dima{The ability to express automorphic forms through Whittaker-Fourier coefficients also implies that minimal and next-to-minimal representations cannot appear in the cuspidal spectrum.} For split classical groups of rank bigger than 2 the latter follows from Theorem \ref{thm:Glob}\eqref{it:GlobCusp} and Corollary \ref{cor:GlobCusp}.  }

\DimaA{
The following table compares our notation to the notation of several papers regarding invariants  of local and global representations.

\begin{tabular}{|c|c|c|c|c|c|}\hline
Present paper &  \cite{MW} & \cite{BV,Ross} & \cite{Ginz} & \cite{JLS} & \cite{GGS} \\\hline
$\WO$ & $\cN_{Wh}$ &  &  & $\WF$ &  \\\hline
$\WS$ & $ $ &  & $\cO_G$ &  &  \\\hline
$\WF$ & maximal elements in $\cN_{\mathrm{tr}}$ &  &  &  &  \\\hline
$\overline{\WF}$ &  & $\mathrm{AS}$ &  &  & $\WF$ \\\hline
\end{tabular}
}

\subsection*{Acknowledgements}

We thank Joseph Bernstein, Yuanqing Cai, David Ginzburg, Raul Gomez, Henrik Gustafsson,   David Kazhdan, Axel Kleinschmidt, Erez Lapid, Baiying Liu, Daniel Persson, Gordan Savin, Eitan Sayag, and David Soudry for fruitful discussions.

D.G. was partially supported by ERC StG grant 637912, and ISF grant 249/17.
S.S. was partially supported by Simons Foundation grant 509766.



\subsection{Notation}\label{subsec:not}$\,$
Let $F$  be either $\R$ or a finite extension of $\Q_p$ and let $\fg$ be a reductive Lie algebra over $F$.
We say that an element $S\in \fg$ is \emph{rational semi-simple} if its adjoint action on $\fg$ is diagonalizable with  eigenvalues in $\Q$.
For a rational semi-simple element $S$ and a rational number $r$ we denote by $\fg_r^S$ the $r$-eigenspace of the adjoint action of $S$ and by $\fg_{\geq r}^S$ the sum $\bigoplus_{r'\geq r}\fg_{r'}^S$. We will also use the notation $(\fg^*)_{ r}^S$ and $(\fg^*)_{\geq r}^S$ for the corresponding grading and filtration of the dual Lie algebra $\fg^*$.  For $X\in \fg$ or $X\in \fg^*$ we denote by $\fg_X$ the centralizer of $X$ in $\fg$, and by $G_X$ the centralizer of $X$ in $G$.

If $(f,h,e)$ is an $\sl_2$-triple in $\fg$, we will say that $e$ is a nil-positive element for $h$, $f$ is a nil-negative element for $h$, and $h$ is a neutral element for $e$. For a representation $V$ of $(f,h,e)$ we denote by $V^e$ the space spanned by the highest-weight vectors and by $V^f$ the space spanned by the lowest-weight vectors. We refer to \cite[\S 11]{Bou} or \cite{KosSl2} for  standard facts on $\sl_2$-triples.
 We will say that $h\in \fg$ is a \emph{neutral} element for $\varphi\in \fg^*$ if $h$ can be completed to an $\sl_2$-triple $(f,h,e)$ such that $\varphi$ is given by the Killing form pairing with $f$. By the Jacobson-Morozov theorem such $h$ exists and is unique up to conjugation by the centralizer of $\varphi$.

From now on we fix a non-trivial unitary additive character $\chi:\F\to \mathrm{S}^1$ such that\\
\begin{equation}\label{=chi}
\begin{split}
&\text{if }\F=\R \text{ we have }\chi(x)=\exp(2\pi i x) \text{ and }\\
&\text{if }\F \text{ is non-Archimedean the kernel of }\chi \text{ is the ring of integers}.
\end{split}
\end{equation}


For non-archimedean $F$ we will work with  $l$-groups, {\it i.e.} Hausdorff topological groups having a basis for the topology at the identity consisting of open compact subgroups. This generality includes $F$-points of algebraic groups defined over $F$, and their finite covers (see \cite{BZ}).

For $F=\R$ we will work with affine Nash groups, {\it i.e.} Lie groups that are given in $\R^n$ by semi-algebraic equations, as well as the graphs of their multiplication maps.  This generality includes $\R$-points of algebraic groups defined over $\R$, and their finite covers (see \cite{dCl,AG,Sun,FS}).

\begin{notn}\label{not:rep}
If $G$ is an $l$-group, we denote by $\Rep^{\infty}(G)$ the category of smooth representations of $G$ in complex vector spaces. For $V,W\in \Rep^{\infty}(G)$, $V\otimes W$ will denote the usual tensor product over $\C$. We denote by $\cM(G)\subset \Rep^{\infty}(G)$ the subcategory consisting of representations of finite length.

If $G$ is an affine Nash group, we denote by $\Rep^{\infty}(G)$ the category of smooth nuclear \Fre representations of $G$ of moderate growth. This is essentially the same definition as in \cite[\S 1.4]{dCl} with the additional assumption that the representation spaces are nuclear (see e.g. \cite[\S 50]{Tre}). For $V,W\in \Rep^{\infty}(G)$, $V\otimes W$ will denote the completed projective tensor product. We denote by $\cM(G)\subset \Rep^{\infty}(G)$ the subcategory consisting of admissible finitely generated representations, see \cite[Ch. 11]{Wal}. This category is abelian, see \cite[Ch. 11]{Wal} or \cite{CasGlob}.
\end{notn}

For $\pi\in \Rep^{\infty}(G)$, denote by $\pi_G$ the space of coinvariants, i.e. quotient of $\pi$ by the intersection of kernels of all $G$-invariant functionals. Explicitly,
$$\pi_G=\pi/\overline{\{\pi(g)v -v\, \vert \,v\in \pi, \, g\in G\}}, $$
where the closure is needed only for Archimedean $F$. In the latter case,
for connected $G$ we have  $\pi_G=\pi/\overline{\fg_{\C}\pi}$ which in turn is equal to the quotient of $\oH_0(\fg,\pi)$ by the closure of zero.

\begin{defn}\label{def:ind}
If $G$ is an $l$-group, $H\subset G$ a closed subgroup and $\pi\in \Rep^{\infty}(H)$, we denote by $\ind_H^G(\pi)$ the smooth compactly-supported induction as in \cite[\S 2.22]{BZ}.

If $G$ is an affine Nash group, $H\subset G$ a closed  Nash subgroup and $\pi\in \Rep^{\infty}(H),$ we denote by $\ind_H^G(\pi)$ the Schwartz induction as in \cite[\S 2]{dCl}.

\end{defn}

%
%
%
%
%
%

This induction has the expected properties of small induction, such as induction by stages and  Frobenius reciprocity.

\section{Degenerate Whittaker models and Fourier coefficients}\label{sec:DegWhit}

\subsection{Definitions}\label{subsec:DefMod}
Let $G$ be a finite central extension of the group $G^{\mathrm{alg}}$ of $F$-points of a reductive algebraic group defined over $F$. Let $G^{ad}$ denote the corresponding adjoint algebraic group.

\begin{lem}[{\cite[Appendix I]{MW_Scr}}]\label{lem:cov}
Let $U\subset G^{\mathrm{alg}}$ be a unipotent subgroup, and $\hat U$ be the preimage of $U$ in $G$. Then there exists a unique open subgroup $U'\subset \hat U$ that projects isomorphically onto $U$.
\end{lem}

We will therefore identify the unipotent subgroups of $G^{\mathrm{alg}}$ with their liftings in $G$.

\begin{defn}
Let $W_n$ denote the  2n-dimensional $\F$-vector space $(\F^n)^* \oplus \F^n$ and let $\omega$ be the standard symplectic form on $W_n$.
The Heisenberg group $H_n$ is the algebraic group with underlying algebraic variety $W_{n} \times \F$ with the group law given by $$(w_1,z_1)(w_2,z_2) = (w_1+w_2,z_1+z_2+1/2\omega(w_1,w_2)).$$
Note that $H_0=\F$.
\end{defn}


\begin{defn}
Let $\chi$ be the additive character  of $\F$, as in \eqref{=chi}.  Extend $\chi$ trivially to a character of the commutative subgroup $0\oplus \F^n\oplus \F \subset H_n$.
The oscillator representation ${\varpi}_\chi$ is the unitary induction of $\chi$ from $0\oplus \F^n\oplus \F$ to $H_n$. Define the smooth oscillator representation
$\sigma_\chi$ to be the space of smooth vectors in ${\varpi}_\chi$.
One can show that $\sigma_{\chi}=\ind_{0\oplus \F^n\oplus \F}^{H_n}(\chi)$.
\end{defn}

\begin{defn}\label{def:DegWhit}
\begin{enumerate}[(i)]
\item A \emph{Whittaker pair} is an ordered pair $(S,\varphi)$ such that $S\in \fg$ is rational semi-simple, and $\varphi\in (\fg^*)^{S}_{-2}$. Given such a Whittaker pair, we define the  \emph{degenerate Whittaker model} $\cW_{S,\varphi}$ in the following way: let $\fu:=\fg_{\geq 1}^S$.
Define an anti-symmetric form $\omega_\varphi$ on $\fg$ by $\omega_\varphi(X,Y):= \varphi([X,Y])$. Let $\fn$ be the radical of $\omega_\varphi|_{\fu}$. Note that $\fu,\fn$ are nilpotent subalgebras of $\fg$, and $[\fu,\fu]\subset \fg^S_{\geq 2}\subset \fn$.
 Let $U:=\Exp(\fu)$ and $N:=\Exp(\fn)$ be the corresponding nilpotent subgroups of $G$.
Let  $\fn' :=\fn\cap \Ker(\varphi), \, N':=\Exp(\fn')$. If $\varphi=0$ we define \begin{equation}\cW_{S,0}:=\ind_{U}^G(\C).\end{equation}
Assume now that $\varphi$ is non-zero.
Then $U/N'$ has a natural structure of a Heisenberg group, and its center is $N/N'$. Let $\chi_{\varphi}$ denote the unitary character of $N/N'$ given by $\chi_{\varphi}(\exp(X)):=\chi(\varphi(X))$.
Let $\sigma_\varphi$ denote the oscillator representation of $U/N'$ with central character $\chi_{\varphi}$, and $\sigma'_\varphi$ denote its trivial lifting to $U$. Define \begin{equation}\cW_{S,\varphi}:=\ind_{U}^G(\sigma'_\varphi).\end{equation}


\item For a nilpotent element $\varphi\in \fg^*$, define the \emph{generalized Whittaker model} $\cW_{\varphi}$  corresponding to $\varphi$ to be $\cW_{S,\varphi}$, where $S$ is a neutral element for $\varphi$ if $\varphi\neq 0$ and $S=0$ if $\varphi=0$. We will also call $\cW_{S,\varphi}$ \emph{neutral degenerate Whittaker model}. Since all neutral elements for $\varphi$ are conjugate by the centralizer of $\varphi$, $\cW_{\varphi}$ depends only on the coadjoint orbit of $\varphi$, and does not depend on the choice of $S$. Thus we will also use the notation $\cW_{\cO}$ for a nilpotent coadjoint orbit $\cO\subset \fg^*$. See \cite[\S 5]{GGS} for a formulation of this definition without choosing $S$.


\item To $\pi\in \Rep^{\infty}(G)$ associate the degenerate and generalized Whittaker quotients  by
\begin{equation}\pi_{S,\varphi}:=(\cW_{S,\varphi}\otimes\pi)_G \text{ and }\pi_{\varphi}:=(\cW_{\varphi}\otimes\pi)_G.
\end{equation}
\end{enumerate}
\end{defn}

Slightly different degenerate Whittaker models are considered in \cite{GGS} and denoted
$\cW_{S,\varphi}(\pi)$. They relate to $\pi_{S,\varphi}$ by $\cW_{S,\varphi}(\pi)=\pi_{S,\varphi}^*$.

\subsection{Fourier coefficients}\label{subsec:FC}

Let $K$ be a number field and let $\A=\A_{K}$ be its ring of adeles. In this section we let $\chi$ be a non-trivial unitary character  of $\A$, which is trivial on $K$. Then $\chi$ defines an isomorphism between $\A$ and $\hat{\A}$ via the map $a\mapsto \chi_{a}$, where $\chi_{a}(b)=\chi(ab)$ for all $b\in \A$. This isomorphism restricts to an  isomorphism
\begin{equation}\label{eq:chi_isomorphism}
 \widehat{\A/K}\cong \{\psi\in \hat{\A}\, | \psi|_{K}\equiv 1\}=\{\chi_{a}\, | \, a\in K\}\cong K.
\end{equation}
Given an algebraic group $\bfG$ defined over $K$ we will denote its Lie algebra by  $\fg$ and we will denote the group of its adelic (resp. $K$-rational) points by $\bfG(\A)$ (resp. $\bfG(K)$). We will also define the Lie algebras $\fg(\A)$ and $\fg(K)$ in a similar way.

Given a Whittaker pair $(S,\varphi)$ on $\g(K)$, we set $\fu=\fg_{\geq 1}^{S}$ and $\fn$ to be the radical of the form $\omega_{\varphi}|_{\fu}$, where $\omega_{\varphi}(X,Y)=\varphi([X,Y])$, as before. Let $\bfU=\exp \fu$, and  $\bfN=\exp \fn$. Define a character  $\chi_{\varphi}$ on $\bfU(\A)$ by  $\chi_{\varphi}(\exp X)=\chi(\varphi(X))$ and note that it is automorphic, \emph{i.e.} trivial on $\bfU(K)$.
Let $G$ be a finite central extension of $\bfG(\A)$, such that the cover $G\onto \bfG(\A)$ splits over $\bfG(K)$. Fix a discrete subgroup $\Gamma\subset G$ that projects isomorphically onto
$\bfG(K)$. Note that $\bfU(\A)$ has a canonical lifting into $G$, see e.g. \cite[Appendix I]{MW_Scr}.

\begin{definition}
 Let $(S,\varphi)$ be a Whittaker pair for  $\g(K)$ and let $\bfU, \bfN,$ and $\chi_{\varphi}$  be as above. For an automorphic \Dima{form} $f$, we define its \emph{$(S,\varphi)$--Fourier coefficient} to be
\begin{equation}\label{eq:Whittaker-Fourier-coefficient}
\DimaA{\cF}_{S,\varphi}(f):=\int_{N(\A)/N(K)}\chi_{\varphi}(n)^{-1}f(n)dn.
\end{equation}
\end{definition}
\Dima{
Observe that $\DimaA{\cF}_{S,\varphi}$ defines a linear functional on the space of automorphic forms. For a  subrepresentation $\pi$ of the space of automorphic forms on $G$,   we  denote the restriction of $\DimaA{\cF}_{S,\varphi}$ to $\pi$ by $\DimaA{\cF}_{S,\varphi}(\pi)$. We denote by $\WO(\pi)$ the set of $\bfG(K)$-orbits of all $\varphi\in \fg(K)$ with $\DimaA{\cF}_{\varphi}(\pi)\neq 0$, and by $\WS(\pi)$ the set of maximal orbits in $\WO(\pi)$. \DimaA{In \cite{GGS}, $\cF_{S,\varphi}$ is denoted $\cW\cF_{S,\varphi}$.}

\begin{remark}
\begin{enumerate}[(i)]
\item The set $\WS(\pi)$ \DimaD{may depend} on the embedding of $\pi$ into the space of automorphic forms.
\item One can define abstract (degenerate) Whittaker quotients for automorphic representations, following Definition \ref{def:DegWhit} \DimaC{and \cite[\S 5.2]{SH}}. The non-vanishing of such a quotient follows from the non-vanishing of the corresponding Fourier coefficient (for some realization), but in general does not imply it. Indeed, \cite{GS} construct an automorphic cuspidal representation of the metaplectic double cover of $\SL_2(\A)$ that has an abstract Whittaker model corresponding to some non-zero nilpotent $\varphi$, but with vanishing $\DimaA{\cF}_{\varphi}$.
One usually considers Fourier coefficients rather than abstract quotients since the former have number-theoretic and string theoretic applications.

\item To the best of our knowledge, the Whittaker support and the non-vanishing of abstract quotients are the only ways to attach nilpotent orbits to automorphic representations. Some authors call the \DimaA{$\WO(\pi)$} \emph{the wave-front set} just by definition.
\end{enumerate}
\end{remark}
}
\subsection{Comparison between different Whittaker pairs}\label{subsec:compar}


%
%

We start with a corollary of the Stone-von-Neumann theorem.
\begin{cor}[{See e.g. \cite[Corollary 2.4.5]{GGS}}]\label{cor:StvN}
Let $W:=(\F^n)^* \oplus \F^n$ and let $H$ be the corresponding Heisenberg group and $\sigma_{\chi}$ its \DimaC{smooth} oscillator representation.
 Let $L\subset W$ be a Lagrangian subspace. Extend $\chi$ trivially from $F$  to the abelian subgroup $L\oplus \F \subset H_n$. Then $$\ind_{L\oplus \F}^{H_n}\chi \cong \sigma_\chi.$$
\end{cor}
Using a version of Frobenius reciprocity, and induction by stages we get:

%

\begin{lem}[{See e.g. \cite[Lemma 2.5.3]{GGS2}}]\label{lem:WhitFrob}
Let $\fl\subset \fg_{\geq 1}^S$ be a maximal isotropic subalgebra and $L:=\Exp(\fl)$. Let $\pi\in \Rep^{\infty}(G)$. Then
$$\cW_{S,\varphi} \cong \ind_{L}^G\chi_{\varphi} \text{ and } \pi_{S,\varphi}\cong (\pi\otimes\chi_{\varphi})_L.$$
\end{lem}

This lemma is the central tool for the comparison of different degenerate Whittaker models.
\DimaA{
\begin{remark}\label{rem:RE}
\begin{enumerate}[(i)]
\item
The previous lemma implies that for any two maximal isotropic subalgebras $\fl,\fr\subset \fg_{\geq 1}^S$ we have an isomorphism
\begin{equation}\label{=Heis}
\ind^{G}_{\Exp(\fl)}\chi_\varphi\simeq \ind^{G}_{\Exp(\fr)}\chi_\varphi.
\end{equation}
\item A similar isomorphism was  introduced in \cite[Lemma 2.2]{GRS} and called \emph{root exchange} (see also \cite[Lemma 7.11]{GRS_Book} or \cite[Lemma A.1]{LapidMao}).

\item Observe that the isomorphism \eqref{=Heis} can be realized explicitly as an integral transform: given $f\in  \ind^{G}_{\Exp(\fl)}\chi_\varphi$ we can define $\check{f}\in  \ind^{G}_{\Exp(\fr)}\chi_\varphi$ simply by setting
\begin{equation}\label{eq:RE}
\check{f}(g)=\int_{\Exp(\fl\cap\fr\cap \Ker \varphi) \backslash \Exp(\fr\cap \Ker \varphi)}f(ng)\, dn =\int_{\Exp(\fl\cap\fr) \backslash \Exp(\fr)}\chi_{\varphi}(n)^{-1} f(ng)\, dn.
\end{equation}
The previous results imply that the map $f\mapsto \check{f}$ defines an isomorphism  $\ind^{G}_{\Exp(\fl)}\chi_\varphi \IsoTo \ind^{G}_{\Exp(\fr)}\chi_\varphi$. The corresponding formula in \cite[Remark 3.2.2]{GGS} has a typo.
\end{enumerate}
\end{remark}
}

\begin{example}\label{ex:GL} Let us now give several examples for $\GL_n(F)$. We identify $\fg$ with $\fg^*$ using the trace form. For $1\leq i,j\leq n$ let $E_{ij}$ denote the corresponding elementary matrix. Let $B$ denote the
group of upper-triangular matrices and $B'=[B,B]$ denote its unipotent radical.
\begin{enumerate}[(i)]

\item  \label{it:GL2} First, let $n:=2$. Then there are two nilpotent orbits: the zero one and the regular one. For a regular nilpotent $\varphi \in \fg^*$, any corresponding degenerate Whittaker model is the classical Whittaker model, for any $G$. For the zero $\varphi$, there are always several choices. For example, we can always choose $S=0$ and have $U=L=\{\Id\}$. We can also choose $S$ to be any rational semi-simple element. For $\GL_2(F)$ we can take $S=\diag(1,-1)$ and get $U=L=B'$. The corresponding Whittaker quotients are different both in the real and in the p-adic case.
Indeed, $\pi_0=\pi$ always, while $\pi_{S,0}$ is finite-dimensional \DimaC{for all admissible $\pi$ of finite length}. For non-Archimedean $F$ and cuspidal $\pi$ it even vanishes.

\item Let us now take $G=\GL_3(F)$, and let $\varphi:=E_{21}$. Then we can take $ h:=\diag(1,-1,0)$ and $S:=\diag(2,0,-2)$. Then $(h,\varphi)$ is a neutral Whittaker pair, and $(S,\varphi)$ is a Whittaker pair. For $S$ we have $L=U=B'$ and for $h$ we can choose

\begin{equation}
\fl= \left(  \begin{array}{ccc}
     0 & * & *  \\
     0 & 0 & 0  \\
     0 & 0 & 0  \\
   \end{array}  \right) \subset \fu= \left(  \begin{array}{ccc}
     0 & * & *  \\
     0 & 0 & 0  \\
     0 & * & 0  \\
   \end{array}  \right)
\end{equation}

\item For $G=\GL_4(F)$, and let $\varphi:=E_{21}+E_{43}$.  Then we can take $ h:=\diag(1,-1,1,-1)$ and $S:=\diag(3,1,-1,-3)$. Then $(h,\varphi)$ is a neutral Whittaker pair, and $(S,\varphi)$ is a Whittaker pair. For $S$ we have $L=U=B'$ and for $h$ we have
$$ \fl=\fu=
\left(
   \begin{array}{cccc}
     0 & * & 0 & * \\
     0 & 0 & 0 & 0 \\
     0 & * & 0 & * \\
     0 & 0 & 0 & 0 \\
   \end{array}
 \right).
$$
\end{enumerate}
\end{example}
Here we see a new phenomenon: this subalgebra does not lie in the Lie algebra of $B'$, and cannot be conjugated into it  by the centralizer of $\varphi$. However, it is still possible to present the corresponding degenerate Whittaker model as a quotient of $\cW_{\varphi}$, as we explain below.

\begin{thm}[{\cite[Theorem A]{GGS} and \cite[Theorem D]{GGS2}}]\label{thm:MaxModAction}
Let $(S,\varphi)$ be a Whittaker pair. Let $\overline{G_S\varphi}$ denote the closure of the orbit of $\varphi$ under the centralizer of $S$. Then
\begin{enumerate}[(i)]
\item For any $\psi\in \overline{G_S\varphi}$ there exists a natural surjection $\nu:\cW_{\psi}\onto \cW_{S,\varphi}$.
\item \label{it:non0} Let $\pi\in \Rep^{\infty}(G)$ such that $G \cdot \varphi \in \WS(\pi)$.  Then $\pi_{S,\varphi}\neq 0$. Moreover, if $F$ is non-Archimedean then the epimorphism $\pi_{\varphi}\onto \pi_{S,\varphi}$ induced by $\nu$  is an isomorphism.
\end{enumerate}
\end{thm}

\begin{remark}
Let us comment on the Archimedean case of \eqref{it:non0}. Already from Example \ref{ex:GL}\eqref{it:GL2} we see that $\pi_{\varphi}$ is not isomorphic to $\pi_{S,\varphi}$ in this case. We conjecture though that they will be isomorphic for unitary $\pi$. Also, in \DimaC{a work in progress we will} show that they will be isomorphic if we replace the usual quotients $\pi_{h,\varphi}$ and $\pi_{S,\varphi}$ by the inverse limits $\lim \limits  _{\ot} \pi/(\fl_{\varphi})^i\pi$, where $\fl$ is a maximal isotropic subalgebra of $\fg_{\geq 1}^h$ (respectively  $\fg_{\geq 1}^h$) and $\fl_{\varphi}$ is the ideal of the universal enveloping algebra generated by $\fl \otimes \varphi$.
\end{remark}



Let us explain the construction of $\nu$ for the case $\psi=\varphi$.
One can show that $S$ can be presented as $h+Z$, where $h$ is a neutral element for $\varphi$ and $Z$ commutes with $h$ and centralizes $\varphi$.
Consider a deformation $S_t=h+tZ$, and denote by $\fu_t$ the sum of eigenspaces of $ad(S_t)$ with eigenvalues at least 1. We call a rational number $0< t <1 $ {\it regular} if $\fu_t=\fu_{t+\eps}$ for any small enough rational $\eps$, and {\it critical} otherwise. Note that there are finitely many critical numbers, and denote them by $t_1<\dots<t_n$. Denote also $t_0:=0$ and $t_{n+1}:=1$.
Choose a Lagrangian $\fm\subset \fg^Z_0\cap \fg^S_{1}$.
 For each $t$ we define two subalgebras $\fl_t,\fr_t\subset \fu_t$ by
 \begin{equation}\label{=lt}
 \fl_t:=\fm+(\fu_t\cap \fg^{Z}_{< 0})+\Ker(\omega|_{\fu_t})\text{ and }\fr_t:=\fm+(\fu_t\cap \fg^{Z}_{> 0})+\Ker(\omega|_{\fu_t}).
\end{equation}


Both $\fl_t$ and $\fr_t$ are maximal isotropic subspaces with respect to the form $\omega_{\varphi}$, and thus  the restrictions of $\varphi$ to $\fl_t$ and $\fr_t$ define characters of these subalgebras.
Let $L_t:=\Exp(\fl_t)$ and $R_t:=\Exp(\fr_t)$ denote the corresponding subgroups and $\chi_\varphi$ denote their characters defined by $\varphi$.
By Lemma \ref{lem:WhitFrob} we have
$$\cW_{S_t,\varphi}\simeq \ind_{L_t}^G(\chi_{\varphi})\simeq \ind_{R_t}^G(\chi_{\varphi}).$$
We show that for any $0\leq i\leq n, \, \fr_{t_i}\subset \fl_{t_{i+1}}$. This gives a natural epimorphism
$$\cW_{S_{t_i},\varphi}\simeq \ind_{L_{t_{i}}}^G(\chi_{\varphi})\onto \ind_{R_{t_i}}^G(\chi_{\varphi})\simeq \cW_{S_{t_{i+1}},\varphi}.$$
Altogether, we get
\begin{equation}\label{=GGS}
\cW_{h,\varphi}=\cW_{S_{t_0},\varphi}\onto \cW_{S_{t_1},\varphi}\onto \cdots \onto \cW_{S_{t_{n+1}},\varphi}=\cW_{S,\varphi}.
\end{equation}

\begin{example}\label{ex:GLSame}
Let $G:=\GL(4,\F)$ and let $S$ be the diagonal matrix $\diag(3,1,-1,-3)$. Identify $\fg$ with $\fg^*$ using the trace form and let
$f:=\varphi:=E_{21}+E_{43}$, where $E_{ij}$ are elementary matrices.
Then we have $S=h+Z$ with $h=\diag(1,-1,1,-1)$ and $Z=\diag(2,2,-2,-2)$.
Thus $S_t=\diag(1+2t,-1+2t,1-2t,-1-2t)$ and the weights of $S_t$ are as follows:
$$\left(
   \begin{array}{cccc}
     0 & 2 & 4t & 4t+2 \\
     -2 & 0 & 4t-2 & 4t \\
     -4t & -4t+2 & 0 & 2 \\
     -4t-2 & -4t & -2 & 0 \\
   \end{array}
 \right).$$
The  critical numbers are $1/4$ and $3/4$. For $t\geq3/4$, the  degenerate Whittaker model $\cW_{S_t,\varphi}$ is the induction $\ind_{B'}^G\chi_{\varphi}$, where $B'$ is the group of upper-unitriangular matrices.
The sequence of inclusions $\fr_{0}\subset \fl_{1/4}\sim \fr_{1/4}\subset \fl_{3/4}=\fr_{3/4}$  is:

\begin{equation}\label{=ThmEx}
\left(
   \begin{array}{cccc}
     0 & - & 0 & - \\
     0 & 0 & 0 & 0 \\
     0 & - & 0 & - \\
     0 & 0 & 0 & 0 \\
   \end{array}
 \right) \subset \left(
   \begin{array}{cccc}
     0 & - & a & - \\
     0 & 0 & 0 & a \\
     0 & * & 0 & - \\
     0 & 0 & 0 & 0 \\
   \end{array}
 \right)\sim
\left(
   \begin{array}{cccc}
     0 & - & * & - \\
     0 & 0 & 0 & * \\
     0 & 0 & 0 & - \\
     0 & 0 & 0 & 0 \\
   \end{array} \right)
\subset
\left(
   \begin{array}{cccc}
     0 & - & - & - \\
     0 & 0 & * & - \\
     0 & 0 & 0 & - \\
     0 & 0 & 0 & 0 \\
    \end{array} \right)
\end{equation}
   Here, both $*$ and $-$ denote arbitrary elements. $-$ denotes the entries in $\fv_t$ and $*$ those in $\fw_t=\fg_{1}^{S_t}$. The letter $a$ denotes an arbitrary element, but the two appearances of $a$ denote the same numbers.
   The  passage from $\fl_{1/4}$ to $\fr_{1/4}$ is denoted by $\sim$. At  $3/4$ we have $\fl_{3/4}=\fr_{3/4}.$

\end{example}

Let us now sketch the proof of part \eqref{it:non0} of Theorem \ref{thm:MaxModAction}.
The sequence of epimorphisms \eqref{=GGS} naturally defines a sequence of epimorphisms
\begin{equation}\label{=EpiSeq}
\pi_{h,\varphi}=\pi_{S_{t_0},\varphi}\onto \pi_{S_{t_1},\varphi}\onto \cdots \onto \pi_{S_{t_{n+1}},\varphi}=\pi_{S,\varphi}.
\end{equation}
We see that for each $i$, $\pi_{S_{t_{i+1}},\varphi}$ is the quotient of $\pi_{S_{t_i},\varphi}$ by the group $A_i:=L_{t_{i+1}}/R_t$, that we show to be commutative. Using Fourier transform on this group, one shows that in order to prove  the theorem it is enough to show that   $\pi_{S_{t_i},\varphi}$ is a non-generic representation of $A_i$. For that purpose we show that every unitary character of $A_i$ is given by some $\varphi'\in\fg^*$ with $ad^*(S_{t_{i+1}})\varphi'=-\varphi'$ such that $\varphi'$ does not lie in the tangent space to $\cO$ at $\varphi$.
We then define a quasi-Whittaker quotient  $\pi_{S_{t_{i+1}},\varphi,\varphi'}$, and show that its dual is the space of $(A_i,\chi_{\varphi'})$-equivariant functionals on $\pi_{S_{t_i},\varphi}$. Then we generalize \eqref{=EpiSeq} to quasi-Whittaker quotients, construct some additional epimorphisms and deduce the vanishing of $\pi_{S_{t_{i+1}},\varphi,\varphi'}$ from the vanishing of $\pi_{\cO'}$ for all $\cO'\neq \cO$ with $\cO\subset \overline{\cO'}$.

%


Let us now follow this argument in the setting of Example \ref{ex:GLSame}. Let $\pi\in \Rep^{\infty}(G)$ with $G\varphi\in \WS(\pi)$.
The sequence of epimorphisms \eqref{=EpiSeq} is given by the sequence of inclusions \eqref{=ThmEx}. To see that these epimorphisms are non-zero (and are isomorphisms for $F\neq \R$) we need to analyze the dual spaces to $(\fg^{S_{1/4}}_1)^f$ and $(\fg^{S_{3/4}}_1)^f$.
These spaces are spanned by $E_{13}+E_{24}$ and by $E_{23}$ respectively. Thus, the dual spaces are spanned   by $E_{31}+E_{42}$ and by $E_{32}$ respectively.
Note that the joint centralizer of $h,Z$ and $\varphi$ in $G$ acts on these spaces by scalar multiplications, identifying all non-trivial elements.
It is enough to show that $\pi_{S_{1/4},\varphi,E_{31}+E_{42}}=0$ and $\pi_{S_{3/4},\varphi,E_{32}}=0$.

First assume by way of contradiction that $\pi_{S_{3/4},\varphi,E_{32}}\neq 0$. Note that $E_{32}\in \fg^{S_{1}}_{-2}$ and that $\fw_1=0 $. Thus $\fu_1=\fl_1=\fr_{3/4}$ and $$\pi_{S_{1},\varphi+E_{32}}\simeq \pi_{S_{3/4},\varphi,E_{32}}\neq 0.$$ Note that $\Phi:=\varphi+E_{32}=E_{21}+E_{43}+E_{32}$ is a regular nilpotent element, and $S_{1}=S=\diag(3,1,-1,-3)$ is a neutral element for it. Thus $\pi_{\Phi}\neq 0,$ contradicting the assumption that $\Phi$ is maximal in $\WS(\pi)$.

Now assume by way of contradiction that $\pi_{S_{1/4},\varphi,E_{31}+E_{42}}\neq 0$. Note that $E_{31}+E_{42}\in \fg^{S_{1/2}}_{-2}$ and that $\fw_{1/2}=0 $. Thus $\fl_{1/2}=\fu_{1/2}=\fr_{1/4}$ and  $$\pi_{S_{1/2},\varphi+E_{31}+E_{42}}\simeq \pi_{S_{1/4},\varphi,E_{31}+E_{42}}\neq 0.$$ Note that $\Psi:=\varphi+E_{31}+E_{42}=E_{21}+E_{43}+E_{31}+E_{42}$ is a regular nilpotent element, and $S_{1/2}=\diag(2,0,0,-2)$ is a neutral element for it. Thus $\pi_{\Psi}\neq 0,$ contradicting the assumption that $G \cdot \varphi$ is maximal in $\WO(\pi)$.

One can deduce from Theorem \ref{thm:MaxModAction} that for non-Archimedean $F$ and quasi-cuspidal $\pi\in \Rep^{\infty}(G)$, the  orbits in $\WS(\pi)$ are $F$-distinguished.


\begin{proof}[Sketch of proof of Theorem  \ref{thm:cuspTemp}\eqref{it:cusp}]
Let $\pi$ be quasi-cuspidal and let $\cO\in \WS(\pi)$. Suppose by way of contradiction that $\cO$ is not $F$-distinguished. Thus there exists a proper parabolic subgroup $P\subset G$, a
Levi subgroup $L\subset P$ and a nilpotent $f\in \fl$ such that $\varphi\in \cO$, where $\varphi\in \fg^*$ is given by the Killing form pairing with $f$.  Let $h$ be a neutral element for $f$ in $\fl$.
Choose a rational-semisimple element $Z\in \fg$ such that $L$ is the centralizer of $Z$, $\fp:=\fg^{Z}_{\geq 0}$ is the Lie algebra of $P$, and all the positive eigenvalues of $Z$ are bigger than all the eigenvalues of $h$ by at least 2. Note that $\fn:=\fg^{Z}_{> 0}$ is the nilradical of $\fp$. Let $S:=h+Z$. By construction we have $\fn\subset \fg^S_{>2}$ and thus the degenerate Whittaker quotient $\pi_{S,\varphi}$ is a quotient of $r_P\pi$. By Theorem \ref{thm:MaxModAction}, the maximality of $\cO$ implies $\pi_{S,\varphi}\simeq \pi_{\varphi}$.
Thus $r_P\pi$ does not vanish, in contradiction with the condition that $\pi$ is quasi-cuspidal.
\end{proof}

In the global case, we have the following analog of Theorem \ref{thm:MaxModAction}.

\begin{thm}[{\cite[Theorem C]{GGS} and \cite[Theorem 8.0.3]{GGS2}}]\label{thm:GlobMaxModAction}
Let $(S,\varphi)\in \fg(K)\times \fg^*(K)$ be a Whittaker pair. Let $\pi$ be an automorphic representation of $G$.  Then
\begin{enumerate}[(i)]
\item \label{it:SmallOrb} The functional $\DimaA{\cF}_{S,\varphi}(\pi)$ can be obtained by a series of integral transforms from $\DimaA{\cF}_{\varphi}(\pi)$.
\item \label{it:Gnon0} If ${\bf G}(K)\varphi \in \WS(\pi)$ then the functional $\DimaA{\cF}_{\varphi}(\pi)$ can be obtained by a series of integral transforms from $\DimaA{\cF}_{S,\varphi}(\pi)$.
\end{enumerate}
\end{thm}

To prove this analog, we argue in the same way, but replace Lemma \ref{lem:WhitFrob} by an explicit integral transform, in the spirit of Remark \ref{rem:RE}.

This theorem can be generalized. Namely, for \eqref{it:Gnon0} we show in \cite{GSPhys1}  that  if $\varphi \notin \WS(\pi)$ then $\DimaA{\cF}_{\varphi}(\pi)$ can be obtained by a series of integral transforms from $\DimaA{\cF}_{S,\varphi}(\pi)$ and from $\{\DimaA{\cF}_{\psi} \, \vert \, {\bf G}(K)\psi \in \WS(\pi)\}$. For \eqref{it:SmallOrb}  let $\gamma=(e,h,f)$ be an $\sl_2$-triple in $\fg$, let $Z$ commute with $\gamma$ and let $f'\in \fg^e\cap \fg_{-2}^{h+Z}$. Let $\varphi,\psi\in \fg^*$ be given by Killing form pairings with $f,f+f'$ respectively. We \DimaC{show in}  \cite{GGS2} that $\DimaA{\cF}_{\varphi}(\pi)$ can be obtained by a series of integral transforms from $\DimaA{\cF}_{K+h,\psi}(\pi)$.
\DimaC{We deduce that for ${\bf G}=\GL_n$, the set  $\WO(\pi)$ is closed under the natural ordering on orbits.}

\subsection{The Slodowy slice}
Let us define an important  geometric notion used in several proofs in \cite{GGS,GGS2}. For an $\sl_2$-triple $(e,h,f)$ in $\fg$, define the Slodowy slice at $f$ to the orbit $\cO$ of $f$ to be the affine space $f+\fg^e$. This space is strongly transversal to $\cO$. Indeed, the tangent space at $f$ to $\cO$ is $[f,\fg]$ and $\fg=[f,\fg]+\fg^e$.

Quantizing the Slodowy slice one gets a finite $W$-algebra (see \emph{e.g.} \cite{Pre,GG}). The representation theory of finite $W$-algebras seems to be closely related to the theory of generalized Whittaker models over Archimedean fields.

\section{Admissible and quasi-admissible orbits}\label{sec:PfAdm}

\subsection{Definitions}\label{subsec:cov}

Let $\gamma=(e,h,f)$ be an $\sl_2$-triple in $\fg$ and let  $\varphi \in \g^{\ast}$ be given by the Killing form pairing with $f$.
 Let $G_{\gamma}$ denote the joint centralizer of the three elements of $\gamma$. It is well known that $G_{\gamma}$ is a Levi subgroup of $G_{\varphi}$.
 Recall that $\varphi$ induces a non-degenerate symplectic form $\omega_{\varphi}$ on $\g^h_1$ and note that $G_{\gamma}$ acts on $\g^h_1$ preserving the symplectic form. That is, there is a natural map $G_{\gamma}\rightarrow \Sp(\g^h_1)=\Sp(\omega_{\varphi})$. Let $\widetilde{\Sp(\omega_{\varphi})}\rightarrow \Sp(\omega_{\varphi})$ be the metaplectic double covering, and set
\[
 \widetilde{G_{\gamma}}=G_{\gamma}\times_{\Sp(\omega_{\varphi})} \widetilde{\Sp(\omega_{\varphi})}.
\]
Observe that the natural map $\widetilde{G_{\gamma}}\rightarrow G_{\gamma}$ defines a double cover of $G_{\gamma}$. We denote by $M_{\gamma}$  the  subgroup of $G_{\gamma}$  generated by the unipotent elements.
Let $\widetilde{M_{\gamma}}$ denote the preimage of  $M_{\gamma}$ under the projection $\widetilde{G_{\gamma}}\to G_{\gamma}$. Note that  different choices of $\gamma$ with the same $f$ lead to conjugate groups $G_{\gamma}$ and $M_{\gamma}$.
One can also define a covering $\widetilde{G_{\varphi}}$ of the group $G_{\varphi}$, using the symplectic form defined by $\varphi$ on $\fg/\fg^{\varphi}$. It is easy to see that this cover splits over the unipotent radical of $G_{\varphi}$, and that the preimage of $G_{\gamma}$ in $\widetilde{G_{\varphi}}$ is isomorphic to $\widetilde{G_{\gamma}}$, see {\emph e.g.} \cite{Nev}.

\begin{defn}
Let $\bf H$ be a linear algebraic group defined over $F$, and
fix an embedding $\bf H \into \GL_n$. Denote by $\fh$ the Lie algebra of ${\bf H}(F)$ and by $H_0$ the open normal subgroup of ${\bf H}(F)$
generated by the image of the exponential map $\fh\to {\bf H}(F)$.
\end{defn}
Note that $H_0$ does not depend on the embedding of $\bf H$ into $\GL_n$. Note also that if $\bf H$ is semi-simple then $H_0={\bf H}(F)$ and if $F=\R$ then $H_0$ is the connected component of ${\bf H}(F)$. For $H'$ a  finite central extension  of $H$, we  define $H'_0$ to be the preimage of $H_0$ under the projection $H'\onto H$.

\begin{definition}[{\cite{Nev}}]
 We say that a nilpotent orbit $\Orb\subset \g^{\ast}$ is \emph{admissible} if for some (equivalently, for any) choice of  $\varphi \in \cO$, the covering $\widetilde{G_{\varphi}}\rightarrow G_{\varphi}$ splits over $(G_{\varphi})_0$.
\end{definition}

As observed in \cite{Nev}, this definition of admissibility is compatible with Duflo's original definition for the Archimedean case, given in  \cite{D}.

\begin{definition}
 We say that a nilpotent orbit $\Orb\subset \g^{\ast}$ is \emph{quasi-admissible} if for some (equivalently, for any) 
 $\varphi \in \cO$, the covering $\widetilde{G_{\varphi}}\rightarrow G_{\varphi}$ admits a finite dimensional \emph{genuine} representation, that is, a finite dimensional representation on which the
  non-trivial element $\eps$ in the preimage of $1\in G_{\varphi}$ acts by $-\Id$. 
\end{definition}


\subsection{On the proof of Theorem \ref{thm:adm}}

%

 Let $\gamma=(f,h,e)$ be an $\sl_2$-triple.
 Let $\varphi\in \fg^*$  be given by the Killing form pairing with $f$.
Let us define the action of $\widetilde{G_{\gamma}}$ on $\pi_{\varphi}$.
We will use the notation of Definition \ref{def:DegWhit}.
Since the adjoint action of $G_{\gamma}$ preserves $\fg^h_1$ and the symplectic form on it, it preserves $U/N'$. Since $\sigma_{\varphi}$ is the unique smooth irreducible representation of $U/N'$ with central character $\chi_{\varphi}$, we have a projective action of $G_{\gamma}$ on $\sigma_{\varphi}$. By \cite{Weil} this action lifts to a genuine representation of $\widetilde{G_{\gamma}}$. This gives raise to an action of $\widetilde{G_{\gamma}}$ on $\cW_{\varphi}$
by $(\tilde g f)(x)=\tilde g( f(g^{-1}x))$. This action  commutes with the action of $G$ and thus defines an action of $\widetilde{G_{\gamma}}$ on $\pi_{\varphi}=(\cW_{\varphi}\otimes \pi)_G.$

The technique described in \S \ref{subsec:compar} implies the following theorem.


\begin{thm}[{\cite[Theorem C]{GGS2}}]\label{thm:MaxFin}
Let $\pi\in \Rep^{\infty}(G)$ and assume that $G \cdot \varphi \in \WS(\pi)$.
Then the action of $\widetilde{M_{\gamma}}$ on \DimaD{the dual space} $\pi_{\varphi}^*$ is locally finite. \DimaD{Moreover, if $F$ is non-Archimedean this action is by $\pm \Id$.}
\end{thm}

Since the action  of $\widetilde{M_{\gamma}}$ on $\pi_{\varphi}$ is genuine, this theorem implies Theorem \ref{thm:adm} after some geometric considerations.

\section{Wave-front sets}\label{sec:WF}

\subsection{Definition}\label{subsec:DefWF}

Let $\pi\in\mathcal{M}(G)$. Let $\chi_{\pi}$ be the character of $\pi$. It is a generalized function on $G$ and it defines a generalized function $\xi_{\pi}$ on a neighborhood of zero in ${\mathfrak{g}_{n}}$,
by restriction to a neighborhood of $1\in G$ and applying logarithm.
In the non-Archimedean case $\xi_{\pi}$ is  a combination of Fourier transforms of $G$-measures of nilpotent coadjoint orbits (\cite{HowGL},\cite[p. 180]{HCWF}).
The measures  extend to $\fg^*$ by \cite{RangaRao}.
In the Archimedean case, the leading term of the asymptotic expansion of $\xi_{\pi}$ near 0 is  equal to such a linear combination \cite[Theorems 1.1 and 4.1]{BV}. For each nilpotent orbit $\cO$ denote by $c_{\cO}(\pi)$ the coefficient of the Fourier transform of the appropriately normalized $G$-invariant measure of $\cO$ in the decomposition of $\xi_{\pi}$.

Define the wave front set $\WF(\pi)$ to be the set of orbits $\cO$ such that  $c_{\cO}(\pi)\neq 0$ and $c_{\cO'}=0$ for every orbit $\cO'$ that includes $\cO$ in its closure.
Denote by $\overline{WF}(\pi)$ the closure of the union of all the orbits in $\WF(\pi)$.

The behaviour of $\WF(\pi)$ under induction is studied in \cite{MW,BarBoz}.
It corresponds to the induction of nilpotent orbits defined in \cite[\S II.3]{Spa}. Namely, let $P\subset G$ be a parabolic subgroup and $L$ be the reductive quotient of $P$. Let $\fl$ and $\fp$ denote the Lie algebras of $L$ and $P$. Let $\cO\subset \fl^*\subset \fp^*$ be a nilpotent orbit and let $\hat \cO$ be its preimage under the restriction map $\fg^*\onto \fp^*$. A $G$-orbit $\cO'\subset \fg$ is said to be induced from $\cO\subset \fl$ if $\cO'$ intersects $\hat \cO$ by an open subset. All the induced orbits are conjugate over the algebraic closure of $F$.

The behaviour of $\WF(\pi)$ under cohomological induction is studied in \cite{BVInd}.

It follows from Theorem \ref{thm:MW} and from \cite{SV} that for any $\cO\in \WF(\pi), \, c_{\cO}$ is a natural number (at least for algebraic $G$).  In the $p$-adic case it is shown in \cite[Corollary 1]{Var} that $c_{\cO}\in \Q$ for all $\cO$. Moreover, for $\GL_n(F)$ $c_{\cO}\in \Z$ for all $\cO$ by \DimaD{\cite[\S 10.3]{DeB} (cf. \cite{HowGL})}. If $\pi \in \cM(\GL_n(F))$ is irreducible, then $c_{\cO}=1$ for all $\cO\in \WF(\pi)$ by \cite[\S II.2]{MW}.
\subsection{On the proof of Theorem \ref{thm:MW}}\label{sec:PfMW}
We give a sketch of the proof, which closely follows \cite{MW,Var}  but also highlights some novel features suggested by \cite{Say} (cf. \cite[\S 5.2]{AGS_Zeig}). This is part of ongoing joint work with E. Sayag, R. Gomez and A. Kemarsky.

Let $(H,\varphi)$ be a Whittaker pair such that there exists a morphism $\nu:F^{\times}\to G$ with $d_1\nu(1)=H$. Any neutral Whittaker pair has this property. Let $\cO:=G\cdot \varphi$. Let $\pi \in \cM(G)$ and assume that
\DimaA{
\begin{equation}\label{=max}
\text{for any }\cO'\in \WF(\pi) \text{ with } \cO \subseteq \overline{\cO'}, \text{ we have } \cO'=\cO.
\end{equation}}
It is enough to show that under this assumption we have $\dim \cW_{H,\varphi}=c_\cO$.

Assume for simplicity that all the eigenvalues of $ad(H)$ are even, and that $\varphi\neq 0$.

The starting point of the proof is the technique of approximation of unipotent subgroups of $G$ by open compact subgroups of $G$, as in Jacquet's proof of exactness of his parabolic reduction functors.
More precisely, let $U\subset G$ and $\chi_{\varphi}$ be the nilpotent subgroup and its character constructed in Definition \ref{def:DegWhit}. Then
\cite{MW,Var} construct, following \cite{HowKir,Rod},  a descending sequence of open compact subgroups $K_n$ and their characters $\chi_n$ such that
$\bigcap K_n$ is trivial, $U \subset \bigcup t^nK_nt^{-n}$, and for each $n$, and each $u\in t^nK_nt^{-n}\cap U$,  $\chi_{\varphi}(u)=\chi_n(t^{-n}ut^n)$.

Let $e_n$ denote the Hecke algebra elements given by integration on $K_n$ versus the product of the Haar measure by $\chi_n$. For $n$ big enough,
$K_n$ lies in the image of the exponential map, and we  consider
 the lifting $\exp^*(e_n)$ to $\fg$ and its Fourier transform $\cF(\exp^*(e_n))$.
 Let $p: \fg^*\setminus \{0\}\to \mathbb{P}(\fg^*)$ denote the natural projection to the projective space.
Let $\mu_{\cO}$ denote the  appropriately normalized $G$-invariant measure on $\cO$.
 The $e_n$ are constructed such that
\begin{enumerate}[(a)]
\item $\cF(\exp^*(e_n))$\text{ is the characteristic function of an open compact subset }$B_n \subset \fg^*$.
\item $\mu_{\cO}(\cO\cap B_n)=1$
\item $p(Ad(t)^nB_n)$ converge to $\{p(\varphi)\}$ as $n\to \infty$.
\end{enumerate}
Let $W_n$ be the image of $\pi(e_n)$ and $W'_n :=\pi(t^n)W_n$.
By these properties, assumption \eqref{=max} and the definitions of the character $\chi_{\pi}$ and the wave-front set we have, for $n$ big enough,
\begin{multline}\label{=dim}
\dim W'_n = \dim W_n =  \tr(\pi(e_n))=\chi_{\pi}(e_n)=\sum_{\cO'}c_{\cO'} \mu_{\cO'}(\cF(\exp^*(e_n))=\\
=\sum_{\cO'}c_{\cO'} \mu_{\cO'}(B_n\cap \cO')=0 + \sum_{\cO'\subset \overline{\cO}}c_{\cO'} \mu_{\cO'}(B_n\cap \cO')=0+c_{\cO}+0=c_\cO.
\end{multline}

\DimaA{
It is now left to prove that for $n$ big enough, $W_n'$ projects isomorphically onto $\pi_{H,\varphi}$.
For this purpose} it is shown that for $n$ big enough,  $\pi(t^{n+1}e_{n+1}t^{-n-1})$ defines an embedding $W'_n \into W'_{n+1}$, that is compatible with the projections to $\pi_{H,\varphi}$, and that these projections are also embeddings.

\DimaA{
\begin{equation}\label{eq:fiberdiag}
\xymatrixcolsep{3pc}
\xymatrix{
 W_n'\ar@{^{(}->}[r]
 \ar@{^{(}->}[rd]  & W_{n+1}'\ar@{^{(}->}[d]\ar@{^{(}->}[r] & \dots\ar@{^{(}->}[ld]\\
  & \pi_{H,\varphi}  &
 }
\end{equation}
}

Since $\dim W'_n=c_\cO$, the images of these embeddings coincide and have dimension $c_\cO$. Since $t^nK_nt^{-n}$ ``approximate" $U$, and $\chi_n$ are compatible with $\chi_{\varphi}$, this image is the whole space $\pi_{H,\varphi}$.

\subsection{Archimedean case}\label{subsec:GS}


\Dima{In this subsection we fix $G$ to be a real reductive group, and consider the relationship between wave-front set, Whittaker support, and
several additional invariants. One of these invariants is the  associated variety of the annihilator of $\pi$, that we call for brevity \emph{annihilator variety} and denote $\cV(\pi)$. It is defined for all $\pi\in \Rep^{\infty}(G)$, as the zero set  in $\fg_{\C}^*$ of the ideal in the symmetric algebra $S(\fg_{\C})$, which is generated by the symbols of the  annihilator ideal of $\pi$ in the universal enveloping algebra $U(\fg_{\C})$.
This invariant is more coarse than $\WF(\pi)$, as the following theorem shows.

\begin{thm}[{\cite[Theorem D]{Ross}, using \cite{BB3,Jos}}]
Let $\pi\in \cM(D)$. Then $\cV(\pi)$ is the Zariski closure of $\overline{\WF(\pi)}$. Moreover, if $\pi$ is irreducible then $\cV(\pi)$ is the closure of a single complex nilpotent orbit.
\end{thm}
In particular, if $G$ is a complex reductive group or $G=\GL_n(\R)$
we have $\overline{\WF}(\pi)= \cV(\pi)\cap \fg^*$.
Thus, in these cases $\WF(\pi)$ consists of a single orbit for all irreducible admissible $\pi$. However, for $G=\SL_2(\R)$ this is not the case. Indeed, this group has two regular real nilpotent orbits. For irreducible principal series representations, $\WF$ is the union of these orbits, while for discrete series $\WF$ consists of one of these orbits.

Let us now discuss
the relation of $\cV(\pi)$ and $\WF(\pi)$ to $\WO(\pi)$.}

\begin{thm}[{\cite[Corollary  4]{Mat}}]\label{thm:Mat}
Let $\pi\in \Rep^{\infty}(G)$, and $\cO\in \WO(\pi)$.  Then $\cO\subset \cV(\pi)$.
\end{thm}

Let us now list several results in the other direction.

\begin{thm}[{\cite{MatENS}}]
Suppose that $G$ is a complex reductive group and let $\pi\in \cM(G)$ have regular infinitesimal character. Let $(S,\varphi)$ be a Whittaker pair such that the orbit $G\varphi$ lies in $\WF(\pi)$, and  intersects the nilradical of the parabolic subgroup defined by $S$ in a dense subset. Then $0<\dim \pi_{S,\varphi}<\infty$.
\end{thm}

The result in \cite{MatENS} includes also the vanishing of the corresponding higher homologies.

\begin{thm}[{\cite[\S 3.3]{GGS}}, based on \cite{GS-Gen}]\label{thm:GS}
Let $G$ be quasi-split and algebraic, and let $\pi \in \cM(G)$. Let $(S,\varphi)$ be a Whittaker pair such that $\varphi$ is a principal nilpotent element of a Levi subalgebra of $\fg$. Suppose $\varphi \in \overline{\WF}(\pi)$. Then  there exists $g\in G_{\C}$ such that $ad(g)$ preserves $\fg$ and $\pi_{ad(g)(S),ad(g)(\varphi)}\neq 0$. Moreover, for complex $G$ and for $G=\GL_n(\R)$  we have $\pi_{S,\varphi}\neq 0$.
\end{thm}

The sets $\WF$ and $\WS$ coincide also for representations induced from finite-dimensional representations of parabolic subgroups, by \cite{BarBoz,GSS}. For further results on the equality of $\WF$ and $\WS$  see \cite{GW,Ya01,Pr91,MatDuke,Ma92}.

\Dima{For $\GL_n(F)$ one can express $\pi_{\cO}$ through the Bernstein-Zelevinsky derivatives and their Archimedean analogs, see \cite{BZ-Induced,AGS1} for the definitions  and \cite[Theorem E]{GGS} for their relation to $\pi_{\cO}$. This  implies the following  theorem, for all local fields $F$ with $char F=0$.
\begin{thm}[{\cite[\S II.2]{MW} and \cite[Corollary G]{GGS}}]
Let $\pi\in \cM(\GL_n(F))$ be an irreducible unitarizable representation, and $\cO\in \WF(\pi).$ Then $\dim \pi_{\cO}=1.$
\end{thm}
}

This type of unique models also plays an important role in a family of Rankin-Selberg integrals that represent tensor product L-functions for classical groups (cf. \cite{CFGK,CFGK2}).

\Dima{
In \cite{HoWF}, $\overline{\WF(\pi)}$ is defined for any continuous representation $\pi$ of any Lie group in a Hilbert space. By \cite[Theorem 1.8]{HoWF} and \cite[Theorem 3.4]{Ross}, this definition extends the one given above for reductive $G$ and $\pi\in \cM(G)$ (via realizing $\pi$ as the space of smooth vectors in a Hilbert space representation). The equality of $\overline{\WF(\pi)}$ to several other analytic invariants is proven in \cite[Theorem B]{Ross}.

Let us add that for irreducible unitary representations of type I classical reductive groups, the Howe rank (\cite{HowRank}) is determined by the maximal among the ranks of the matrices lying in the annihilator variety, by \cite{He}.

One can also attach to any smooth admissible representation $\pi$ of a real reductive group $G$ and to any maximal compact subgroup $K\subset G$ an algebraic invariant: \emph{the  associated variety} of the module  $\pi^{(K)}$ of $K$-finite vectors. This variety is a set of nilpotent $K$-orbits in $\fk^{\bot}\subset \fg^*$, with multiplicities, see  \cite{Vo91,Vo17}.
The dimension of this variety equals the Gelfand-Kirillov dimension of $\pi$ and is also equal to half of the dimension of the annihilator variety of $\pi$.

 By  \cite{SV}, if $G$ is algebraic then the maximal orbits in this set correspond to $\WF$ under the Kostant-Sekiguchi correspondence, and the multiplicities are equal to the coefficients $c_{\cO}$. The Kostant-Sekiguchi correspondence maps $K$-orbits in $\fk^{\bot}$ to $G$-orbits in $\fg^*$ that lie in the same orbit of the complexification $G_{\C}$ on $\fg_{\C}^*$, preserving the closure order.}


\begin{thebibliography}{999999}
\bibitem[AG08]{AG} A. Aizenbud, D.
Gourevitch: {\it Schwartz functions on Nash Manifolds,}
International Mathematics Research Notices, Vol. 2008, n.5,
Article ID rnm155, 37 pages. DOI: 10.1093/imrn/rnm155.

\bibitem[AGKLP]{AGKLP} O. Ahlen, H. P. A. Gustafsson, A. Kleinschmidt, B. Liu, D. Persson: \emph{Fourier coefficients attached to small automorphic representations of $SL_n(\A)$,} arXiv:1707.08937.


\bibitem[AGSay15]{AGS_Zeig} A. Aizenbud, D. Gourevitch, E.Sayag:
{\it $\mathfrak{z}$-finite distributions on p-adic groups}, Advances in Mathematics
 \textbf{285},  pp. 1376-1414 (2015).

\bibitem[AGSah15a]{AGS1} A. Aizenbud,  D. Gourevitch, S. Sahi: {\it
Derivatives for smooth representations of $GL(n,{\mathbb{R}})$ and $GL
(n,{\mathbb{C}})$}, Israel J. Math. \textbf{206}, no. 1, 1-38 (2015). See also arXiv:1109.4374.


\bibitem[AGSah15b]{AGS2} A. Aizenbud,  D. Gourevitch, S. Sahi: {\it
Twisted homology of the mirabolic nilradical}, Israel J. Math. \textbf{206}, no. 1, 39-88 (2015). See also arXiv:1210.5389.

\bibitem[BB99]{BarBoz} D. Barbasch, M. Bozicevic: {\it
The associated variety of an induced representation,}
Proc. Amer. Math. Soc. {\bf 127}, no. 1, 279-288 (1999).



%

\bibitem[BB85]{BB3}
W. Borho and J.-L. Brylinski:{\it Differential operators on homogeneous spaces
III}, Inventiones Math. {\bf 80}, 1-68 (1985).

 \bibitem[BV80]{BV}D. Barbasch,  D. A. Vogan: \textit{The local structure of characters. J. Funct. Anal.} \textbf{37}, 27-55 (1980).

\bibitem[BV82]{BVClass}D. Barbasch,  D. A. Vogan: \textit{Primitive ideals and orbital integrals in complex classical groups,} Math. Ann. \textbf{259}, 153-199 (1982).

\bibitem[BV83a]{BVExc}D. Barbasch,  D. A. Vogan: \textit{Primitive ideals and orbital integrals in complex exceptional groups,} J. Alg. \textbf{80}, 350-382 (1983).

\bibitem[BV83b]{BVInd}
D. Barbasch and D. Vogan: \textit{Weyl Group Representations and Nilpotent Orbits}, 21-32, Representation
Theory of Reductive Groups (P.C. Trombi, eds.), Birkhauser, Boston  (1983).

 \bibitem[Bou75]{Bou} N. Bourbaki: \textit{Groupes et algebres de Lie}, Chap. 7 et 8. fasc. XXXVIII. Paris; Hermann (1975).


\bibitem[BZ76]{BZ}I.N. Bernstein, A.V. Zelevinsky: \textit{Representations Of The Group GL(N,F),
Where F Is A Non-Archimedean Local Field,} Uspekhi Mat. Nauk 31:3 , 5-70 (1976).

\bibitem[BZ77]{BZ-Induced}I.N. Bernstein, A.V. Zelevinsky: \textit{Induced
representations of reductive p-adic groups. I.} Ann. Sci. Ec. Norm. Super,
4$^{\text{e}}$serie \textbf{10}, 441-472 (1977).



\bibitem[Car85]{Car}R.W. Carter: \textit{Finite groups of Lie type. Conjugacy
classes and complex characters.} Pure and Applied Mathematics (New York). A
Wiley-Interscience Publication. John Wiley \& Sons, Inc., New York, (1985).
xii+544 pp. ISBN: 0-471-90554-2.

 \bibitem[Cas89]{CasGlob} W. Casselman:  \textit{Canonical extensions of Harish-Chandra modules to representations of G,}
 Can. J. Math., Vol. XLI, No. 3, pp. 385-438 (1989).

\bibitem[CFGK]{CFGK} Y. Cai, S. Friedberg, D. Ginzburg, E. Kaplan
\textit{Doubling Constructions for Covering Groups and Tensor Product L-Functions,}
arXiv:1601.08240.

\bibitem[CFGK2]{CFGK2} Y. Cai, S. Friedberg, D. Ginzburg, E. Kaplan
\textit{Doubling Constructions and Tensor Product L-Functions: the linear case,} arXiv:1710.00905 .


\bibitem[CHM00]{CHM}W. Casselman, H. Hecht, D. Mili\v ci\'c:
\textit{Bruhat filtrations and Whittaker vectors for real groups}. The
mathematical legacy of Harish-Chandra (Baltimore, MD, 1998), 151-190, Proc.
Sympos. Pure Math., \textbf{68}, Amer. Math. Soc., Providence, RI (2000).

\bibitem[CoMG93]{CM}D. Collingwood, W. McGovern: \textit{Nilpotent orbits in
semisimple Lie algebras.} Van Nostrand Reinhold Mathematics Series. Van
Nostrand Reinhold Co., New York, xiv+186 pp (1993).


\bibitem[DeB04]{DeB} S. DeBacker: \textit{Lectures on harmonic analysis for reductive p-adic groups.} In: Representations of Real and
p-adic Groups, Lect. Notes Ser. Inst. Math. Sci. Natl. Univ. Singap., vol. 2, pp. 47-94. Singapore Univ.
Press, Singapore (2004). 

\bibitem[dCl91]{dCl} F. du Cloux: \emph{Sur les repr\'esentations diff\'erentiables des groupes de Lie alg\'ebriques,}  Ann. Sci. \'Ecole Norm. Sup. (4) \textbf{24}, no. 3, 257--318 (1991).


%
%


\bibitem[Duf80]{D} M. Duflo: {\it Construction de repreesentations unitaires d'un groupe de Lie}, in Harmonic
Analysis and Group Representations pp. 129-222, C.I.M.E. (1980).

\bibitem[FS16]{FS}  Y. Fang, B. Sun:
\emph{Chevalley's theorem for affine Nash groups,}
J. Lie Theory \textbf{26}, no. 2, 359-369 (2016).

\bibitem[FGKP18]{FGKP}
P. Fleig, H. P. A. Gustafsson, A. Kleinschmidt, D. Persson: {\it Eisenstein series and automorphic representations,} arXiv:1511.04265,  Cambridge Studies in Advanced Mathematics vol. 176 (2018).

\bibitem[GG02]{GG}
W.L. Gan, V. Ginzburg: {\it
Quantization of Slodowy slices,}
Int. Math. Res. Not., no. 5, 243-255 (2002).

\bibitem[GGJ02]{GGJ} W. T. Gan, N. Gurevich, D. Jiang: {\it Cubic unipotent Arthur parameters and
multiplicities of square integrable automorphic forms,}
Invent. math. {\bf 149}, 225-265 (2002)
DOI 10.1007/s002220200210.




\bibitem[GGS17]{GGS} R. Gomez, D. Gourevitch, S. Sahi: {\it Generalized and degenerate Whittaker models,} Compositio Mathematica \textbf{153}, n. 2, pp. 223-256 (2017).
\href{http://doi.org/10.1112/S0010437X16007788}{DOI:10.1112/S0010437X16007788}.

\bibitem[GGS]{GGS2} R.Gomez, D. Gourevitch, S. Sahi: {\it Whittaker supports for representations of reductive groups,} arXiv:1610.00284v4. To appear in Annales de l'Institut Fourier.

\bibitem[GGKPS1]{GSPhys1} D. Gourevitch, H. Gustafsson, A. Kleinschmidt, D. Persson, S. Sahi:
{\it A reduction principle for Fourier coefficients of automorphic forms,} arXiv:1811.05966.

\bibitem[GGKPS2]{GSPhys2} D. Gourevitch, H. Gustafsson, A. Kleinschmidt, D. Persson, S. Sahi:
{\it Fourier coefficients of minimal and next-to-minimal automorphic representations of simply-laced groups,} arXiv:1908.08296.

\bibitem[Gin06]{Ginz} D. Ginzburg: {\it Certain conjectures relating unipotent orbits to automorphic
representations,} Israel J. Math. {\bf 151}, 323-355 (2006).

\bibitem[Gin14]{Ginz2} D. Ginzburg: {\it Towards
a classification of global integral constructions
and functorial liftings using the small
representations method},
Advances in Mathematics {\bf 254}, 157-186 (2014).




\bibitem[GK75]{GK}
I.~M. Gelfand and D.~A. Kajdan, \emph{Representations of the group
  {${\rm GL}(n,K)$} where {$K$} is a local field}, Lie groups and their
  representations (Proc. Summer School, Bolyai J\'anos Math. Soc., Budapest,
  1971), Halsted, New York, 1975, pp.~95--118.


\bibitem[GMV15]{GMV}
  M.~B.~Green, S.~D.~Miller and P.~Vanhove,
  ``Small representations, string instantons, and Fourier modes of Eisenstein series,''
  J.\ Number Theor.\  {\bf 146} 187-309 (2015)
  doi:10.1016/j.jnt.2013.05.018.





\bibitem[GRS99]{GRS} D. Ginzburg, S. Rallis, D. Soudry: {\it On a correspondence between cuspidal representations of $GL_{2n}$ and $\widetilde{Sp}_{2n}$}. J. Amer. Math. Soc. \textbf{12}, no. 3, 849-907 (1999).

\bibitem[GRS03]{GRS_Sp} D. Ginzburg, S. Rallis, D. Soudry: {\it
On Fourier coefficients of automorphic forms
of symplectic groups,}
Manuscripta math. {\bf 111}, 1-16 (2003).

\bibitem[GRS11]{GRS_Book} D. Ginzburg, S. Rallis, D. Soudry: {\it The descent map from automorphic representations of GL(n) to classical groups,} World Scientific Publishing Co. Pte. Ltd., Hackensack, NJ (2011). x+339 pp.
\bibitem[GS15]{GS-Gen}D. Gourevitch, S. Sahi: \textit{Degenerate Whittaker models for real reductive groups}, American Journal of Mathematics \textbf{137}, n. 2 439-472 (2015).

\bibitem[GS87]{GS} S. Gelbart and D. Soudry: {\it On Whittaker models and the vanishing of Fourier coefficients of cusp forms}, Proc. Indian Acad. Sci. Math. Sci. {\bf 97} , no. 1-3, 67-74 (1987).


\bibitem[GSS18]{GSS}  D. Gourevitch, S. Sahi, E. Sayag: {\it
Analytic continuation of equivariant distributions,} International Mathematics Research Notices, , rnx326 (2018).



\bibitem[GW80]{GW}R. Goodman and N. R. Wallach: \textit{Whittaker vectors and
conical vectors}, J. Funct. Anal. \textbf{39}, no. 2, pp. 199-279 (1980).


\bibitem[GZ14]{GZ} R. Gomez, C.-B. Zhu:
\textit{Local theta lifting of generalized Whittaker models associated to nilpotent orbits},  Geom. Funct. Anal. \textbf{24} , no. 3, 796-853 (2014).



\bibitem[Har12]{Harris}
B. Harris:
{\it Tempered representations and nilpotent orbits,}
Represent. Theory \textbf{16}, 610-619 (2012).


\bibitem[HC77]{HCWF} Harish-Chandra:
\textit{The characters of reductive $p$-adic groups. Contributions to algebra (collection of papers dedicated to Ellis Kolchin)}, pp. 175-182. Academic Press, New York (1977).

\bibitem[He08]{He} H. He: {\it Associated Varieties and Howe's N-Spectrum,}
Pacific Journal Of Math. {\bf 237}, no. 1 (2008).

\bibitem[How74]{HowGL} R. Howe: \textit{The Fourier transform and germs of
characters (case of $GL_{n}$ over a p-adic field)} . Math. Ann 208, 305-322 (1974).

\bibitem[How77]{HowKir} R. Howe: \textit{
Kirillov theory for compact p-adic groups,}
Pacific J. Math. {\bf 73}, no. 2, 365-381 (1977).

\bibitem[How81]{HoWF}
R. Howe: \emph{Wave front sets of representations of Lie groups}, in Automorphic forms, representation theory and arithmetic (Bombay, 1979), pp. 117--140, Tata Inst. Fund. Res. Studies in Math., 10, Tata Inst. Fundamental Res., Bombay, (1981).

\bibitem[How82]{HowRank}R. Howe: \textit{On a notion of rank for unitary
representations of the classical groups}, Harmonic analysis and group
representations pp. 223-331 (1982).


\bibitem[HS16]{SH} J. Hundley, E. Sayag: \emph{Descent Construction
for GSpin Groups},
Memoirs
of the
American Mathematical Society
\textbf{243} n. 1 (2016).






\bibitem[Jia07]{Ji07} D. Jiang: \emph{Periods of automorphic forms}, Proceedings of the International Conference on Complex Geometry and Related Fields,
Studies in Advanced Mathematics \textbf{ 39},   125--148, American Mathematical Society and International Press (2007).

%

\bibitem[JLS16]{JLS} D. Jiang, B. Liu and G. Savin: \emph{Raising nilpotent orbits in wave-front sets}, Represent. Theory {\bf 20}, 419-450 (2016).

\bibitem[Jos85]{Jos}A. Joseph: \textit{On the associated variety of a
primitive ideal}, Journal of Algebra \textbf{93}, n. 2, 509--523 (1985).



\bibitem[Kaz77]{Kaz} D. Kazhdan: {\it Some applications of the Weil representation,} J. Anal. Math. \textbf{32}, 235-248 (1977).

\bibitem[Kaw85]{Ka85} N. Kawanaka: ``Generalized Gelfand-Graev representations and Ennola duality'', in \emph{Algebraic groups and related topics,} Advanced Studies in Pure Math. 6,  pp. 175--206 (1985).

\bibitem[Kos59]{KosSl2}B. Kostant: \textit{The Principal Three-Dimensional Subgroup and the Betti Numbers of a Complex Simple Lie
Group}, American Journal of Mathematics. \textbf{81}, No. 4 973-1032 (1959).

\bibitem[Kos78]{Kos}B. Kostant: \textit{On Whittaker vectors and representation theory}, Invent. Math. \textbf{48}, 101-184 (1978).

\bibitem[LaMao14]{LapidMao} E. Lapid, Z. Mao: \emph{Model transition for representations of metaplectic type}, Int Math Res Notices, doi:10.1093/imrn/rnu225 (2014). See also arXiv:1403.6787.


\bibitem[LoMa15]{LM} H.Y. Loke, J.-J. Ma: \emph{Invariants and $K$-spectrums of local theta lifts}, Compositio Mathematica \textbf{151}, n. 1,  179-206 (2015).


\bibitem[LS08]{LS} H. Y. Loke, G. Savin: {\it On minimal representations of Chevalley groups
of type $D_n$, $E_n$ and $G_2$}. Math. Ann. {\bf 340} (2008), no. 1, 195-208.


\bibitem[Lus79]{LusSpec} G. Lusztig: {\it
A class of irreducible representations of a Weyl group,}
Nederl. Akad. Wetensch. Indag. Math. {\bf 41} , no. 3, 323-335 (1979).

\bibitem[Lus92]{LusFin}  G. Lusztig: {\it
A unipotent support for irreducible representations,}
Adv. in Math. {\bf 94}, 139-179  (1992).


\bibitem[Lyn79]{L} T. E. Lynch: {\it Generalized Whittaker vectors and representation theory}, Thesis, MIT (1979).

\bibitem[Mat87]{Mat}H. Matumoto: \textit{Whittaker vectors and associated
varieties}, Invent. math. \textbf{89}, 219-224 (1987).

\bibitem[Mat90a]{MatENS} H. Matumoto: \emph {$C^{-\infty}$-Whittaker vectors for complex semisimple Lie groups, wave front sets, and Goldie rank polynomial representations}, Ann. Sci. Ecole Norm. Sup. (4) {\bf 23}, no. 2, 311--367 (1990).

\bibitem[Mat90b]{MatDuke}H. Matumoto: \textit{Whittaker modules associated
with highest weight modules}, Duke Math. J. \textbf{60} (1990), no. 1, 59--113.



\bibitem[Mat92]{Ma92} H. Matumoto: \emph {$C^{-\infty}$-Whittaker vectors corresponding to a principal nilpotent orbit of a real reductive linear Lie group, and wave front sets,} Compositio Math. \textbf{82},  189--244 (1992).


\bibitem[MS12]{MilSah} S. D. Miller, S. Sahi: \emph{Fourier coefficients of automorphic forms, character variety orbits, and small representations,} J. Number Theory {\bf 132}, no. 12, 3070--3108 (2012).


\bibitem[Moe96]{Mog}C. Moeglin: \textit{Front d'onde des representations des groupes classiques p-adiques},
Amer. J. Math. \textbf{118}, no. 6, 1313-1346 (1996).

\bibitem[MW87]{MW}C. Moeglin, J.L. Waldspurger: \textit{Modeles de Whittaker
degeneres pour des groupes p-adiques,} Math. Z. \textbf{196}, no. 3, pp 427-452 (1987).

\bibitem[MW95]{MW_Scr}C. Moeglin, J.L. Waldspurger: \textit{Spectral decomposition and Eisenstein series. Une paraphrase de l'Ecriture [A paraphrase of Scripture].} Cambridge Tracts in Mathematics, \textbf{113}. Cambridge University Press, Cambridge (1995). xxviii+338 pp. ISBN: 0-521-41893-3

\bibitem[MT16]{MT}
K. Maktouf, P. Torasso: {\it Restriction de la representation de Weil a un sous-groupe compact maximal,}
J. Math. Soc. Japan {\bf 68}, No. 1, 245-293 (2016).
doi: 10.2969/jmsj/06810245.



\bibitem[Nev99]{Nev} M. Nevins: \textit{Admissible nilpotent coadjoint orbits
of $p$-adic reductive Lie groups}, Representation Theory
{\bf 3}, 105-126 (1999).

\bibitem[Nev02]{Nev2} M. Nevins: \textit{Admissible nilpotent orbits of real
and $p$-adic split exceptional groups},
Representation theory {\bf 6}, 160-189 (2002).

%


\bibitem[NPS73]{NPS73} M. E. Novodvorskii, I. Piatetski-Shapiro: \emph{Generalized Bessel models for a symplectic group of rank $2$}, (Russian) Mat.
Sb. (N.S.)  90 (132),  246--256 (1973).

\bibitem[Oht91a]{Oht} T. Ohta: \textit{Classification of admissible nilpotent orbits in the classical real Lie algebras},
J. Algebra {\bf 136}, no. 2, 290--333 (1991).


\bibitem[PrSk02]{Pre} A. Premet: \textit{Special transverse slices and their enveloping algebras,}
With an appendix by Serge Skryabin.
Adv. Math. {\bf 170} , no. 1, 1-55 (2002).


\bibitem[PT04]{PT} V. L. Popov, E. A. Tevelev:
\textit{Self-dual Projective Algebraic Varieties
Associated With Symmetric Spaces,} Algebraic Transformation Groups and Algebraic Varieties,
Encyclopedia of Mathematical Sciences {\bf  132},
Subseries Invariant Theory and Algebraic Transformation Groups, Vol. III,
Springer-Verlag (2004)


%

 \bibitem[Prz91]{Pr91} T. Przebinda: \emph{Characters, dual pairs, and unipotent representations}, J. Funct. Anal. {\bf 98},  no. 1, 59--96 (1991).

\bibitem[Prz]{Prz} T. Przebinda: \textit{The character and the wave front set correspondence in the stable range,} arXiv:1602.08401.



\bibitem[Ros95]{Ross}W. Rossmann:  {\it Picard-Lefschetz theory and characters of a semisimple Lie group}. Inventiones Mathematicae 121, 579-611 (1995).

\bibitem[RR72]{RangaRao} R. Ranga Rao: {\it
Orbital Integrals in Reductive Groups,} Annals of Mathematics, {\bf 96}, no. 3, 505-510 (1972).

\bibitem[Rod75]{Rod}F. Rodier: \textit{Mod\`{e}le de Whittaker et
caract\`{e}res de repr\'{e}sentations} (French). Non-commutative harmonic
analysis (Actes Colloq., Marseille-Luminy, 1974), pp. 151-171. Lecture Notes
in Math. {\bf 466}, Springer, Berlin, (1975).

\bibitem[Say02]{Say} E. Sayag: {\it A Generalization of Harish-Chandra Regularity Theorem}, Thesis submitted for
the degree of "Doctor of Philosophy", Tel-Aviv University (2002).



\bibitem[Sha74]{Sh74} J. A. Shalika: \emph{The multiplicity one theorem for $\GL_{n}$,} Ann. of Math. 100,  171--193 (1974).




 \bibitem[SV00]{SV}W. Schmid, K. Vilonen: \textit{Characteristic cycles and
 wave front cycles of representations of reductive Lie groups}. Annals of
Mathematics, \textbf{151}, 1071-1118 (2000).

\bibitem[Shen16]{Shen} X. Shen: \textit{Top Fourier Coefficients for Cuspidal Representations of Symplectic Groups,}
International Mathematics Research Notices,
doi: 10.1093/imrn/rnw073 (2016).

\bibitem[Spa82]{Spa} N. Spaltenstein: \textit{Classes unipotentes et sous-groupes de Borel}, Lecture Notes in Math., {\bf 946}, Springer, Berlin-Heidelberg (1982).


\bibitem[Sun15]{Sun} B. Sun: \textit{
Almost linear Nash groups},
Chin. Ann. Math. Ser. B \textbf{36} , no. 3, 355-400 (2015).
\bibitem[Tre67]{Tre} F. Treves: {\it Topological vector spaces, distributions and kernels}, Academic Press, New York-London  xvi+624 pp (1967).

\bibitem[Var14]{Var} S. Varma: {\it On a result of Moeglin and Waldspurger in residual characteristic 2}, Math. Z. \textbf{277}, no. 3-4, 1027--1048 (2014).




\bibitem[Vog91]{Vo91} D. A. Vogan: \textit{Associated varieties and unipotent representations} in \emph{Harmonic analysis on reductive groups} (Brunswick, ME, 1989), pp. 315--388, Progr. Math. {\bf 101}, Birkh\"auser, Boston,  MA (1991).



\bibitem[Vog94]{VogG2}D. A. Vogan: \textit{The Unitary Dual of $G_2$}, Invent. Math. {\bf 116}, 677-791 (1994).

\bibitem[Vog17]{Vo17} D. A. Vogan: \textit{The size of infinite-dimensional representations,} Japan. J. Math. {\bf 12}, 175-210 (2017).

\bibitem[Wall88]{WaJI} N. R. Wallach: \emph{Lie Algebra Cohomology and Holomorphic Continuation of Generalized Jacquet Integrals}, Advanced Studies in Math, vol 14,  123--151 (1988).


\bibitem[Wall92]{Wal}N. Wallach: \textit{Real Reductive groups I,II}, Pure and
Applied Math. \textbf{132}, Academic Press, Boston, MA (1988,1992).

\bibitem[Wei64]{Weil} A. Weil: \textit{Sur certain group d'operateurs unitaires}, Acta
Math. 111,  143--211 (1964).


\bibitem[Yam86]{Ya86} H. Yamashita: \emph{On Whittaker vectors for generalized Gelfand-Graev representations of semisimple Lie groups,} J. Math. Kyoto Univ. 26,  no. 2, 263--298 (1986).

\bibitem[Yam01]{Ya01} H. Yamashita: ``Cayley transform and generalized Whittaker models for irreducible highest weight modules'' in \emph{Nilpotent Orbits, Associated Cycles and Whittaker Models for Highest weight Representations,} Ast\'erisque 273,  81--137 (2001).


\bibitem[Zel80]{Zel}A.V. Zelevinsky : \textit{ Induced representations of
reductive p-adic groups. II. On irreducible representations of Gl(n)}. Ann.
Sci. Ec. Norm. Super, $4^{e}$serie \textbf{13}, 165-210 (1980)


\end{thebibliography}
\end{document}